\documentclass{amsart}
\usepackage{amssymb,amsmath}
\usepackage{multirow}

\newtheorem{tm}{Theorem}
\newtheorem{lm}[tm]{Lemma}
\newtheorem{cor}[tm]{Corollary}
\newtheorem{prop}[tm]{Proposition}

\newtheorem{ex}[tm]{Example}
\newtheorem{rem}[tm]{Remark}
\newtheorem{nots}[tm]{Notations}
\newcommand{\beq}{\begin{equation}}
\newcommand{\eeq}{\end{equation}}
\newcommand{\bga}{\begin{gather*}}
\newcommand{\ega}{\end{gather*}}
\newcommand{\bal}{\begin{align*}}
\newcommand{\eal}{\end{align*}}
\newcommand{\bit}{\begin{itemize}}
\newcommand{\eit}{\end{itemize}}
\newcommand{\btm}{\begin{tm}}
\newcommand{\etm}{\end{tm}}
\newcommand{\blm}{\begin{lm}}
\newcommand{\elm}{\end{lm}}
\newcommand{\bcor}{\begin{cor}}
\newcommand{\ecor}{\end{cor}}
\newcommand{\bprop}{\begin{prop}}
\newcommand{\eprop}{\end{prop}}
\newcommand{\bex}{\begin{ex}}
\newcommand{\eex}{\end{ex}}
\newcommand{\bpr}{\begin{proof}}
\newcommand{\epr}{\end{proof}}
\newcommand{\brem}{\begin{rem}}
\newcommand{\erem}{\end{rem}}
\newcommand{\bdf}{\begin{definition}}
\newcommand{\edf}{\end{definition}}
\newcommand{\bnots}{\begin{nots}}
\newcommand{\enots}{\end{nots}}

\def\C{\mathbb{C}}
\def\N{\mathbb{N}}
\def\R{\mathbb{R}}

\def\Z{\mathbb{Z}}

\def\e{\varepsilon}
\let\a\alpha
\def\trans#1#2{{}^{#1}\mkern-3mu #2}

\def\esssup{\mathop{\rm ess\,sup}}
\def\essinf{\mathop{\rm ess\,inf}}

\def \le {\leqslant}
\def \ge {\geqslant}
\let \<\langle
\let \>\rangle
\let\phi\varphi
\let\kappa\varkappa

\let\o\omega

\def\lpo{L^p(G,\o)}
\def\bSigma{\mathbf{\Sigma}}

\def\supp{\operatorname{supp}}
\begin{document}

\title{Density of translates in weighted $L^p$ spaces on locally compact groups}

\author{Evgeny Abakumov}
\address{Evgeny Abakumov: University Paris-Est, 5 boulevard Descartes, 77454 Marne-la-Vall\'ee, France}
\email{evgueni.abakoumov@u-pem.fr}
\author{Yulia Kuznetsova}
\address{Yulia Kuznetsova: University of Bourgogne Franche-Comt\'e, 16 route de Gray, 25030 Besan\c con, France}
\email{yulia.kuznetsova@univ-fcomte.fr}
\subjclass[2010]{47A16; 37C85; 43A15}
\keywords{locally compact groups; weighted spaces; hypercyclicity; translation semigroups}
\maketitle

\setcounter{tocdepth}{1}

\begin{abstract}
Let $G$ be a locally compact group, and let $1\le p < \infty$. Consider the weighted $L^p$-space
$L^p(G,\omega)=\{f:\int|f\omega|^p<\infty\}$, where $\omega:G\to \R$ is a
positive measurable function. Under appropriate conditions on $\omega$, $G$ acts on $L^p(G,\omega)$
by translations. When is this action hypercyclic, that is, there is a function in this space such that
the set of all its translations is dense in $L^p(G,\omega)$? H.~Salas (1995) gave a criterion of
hypercyclicity in the case $G=\Z$ . Under mild assumptions, we present a corresponding
characterization for a general locally compact group $G$. Our results are obtained in a more general setting when the translations only by a subset $S\subset G$ are considered.
\end{abstract}

\section{Introduction}

Let $G$ be a locally compact group with identity
$e$
and a left Haar
measure $\mu$. Fix  $p, \,1\le p<\infty$, and a weight $\o$, that is, a measurable strictly positive function on $G$ supposed to be  locally $p$-summable.

Consider the weighted space  $L^p(G,\o)=\{f: \int_G |f\o|^pd\mu < \infty\}$ of complex-valued functions on $G$. In this paper the following question is addressed.
Can {\it all} functions in $\lpo$ be approximated arbitrarily well by the left translates of some {\it  single} function $f\in \lpo$?
The same question can be asked if we allow to translate $f$ only by elements of a fixed subset $S\subset G$.

To be more precise, let $S$ be a subset of $G$, and suppose that, for any $s\in S$, the (left)  translation operator

$$(T_s f)(t) = f(s^{-1}t), \,t\in G$$
 is continuous from $L^p(G,\o)$ into itself.
A function $f\in \lpo$ is called {\it $S$-dense} if its $S$-orbit $Orb_S(f) = \{T_s f: s\in S\}$ is dense in $\lpo$. So, we are interested in the existence of $S$-dense vectors in the space $\lpo$.

$S$-density is a particular case of a more general situation, known as the universality phenomenon, see a survey of K.-G.~Grosse-Erdmann \cite{Grosse}.
One of the first examples in this area was given in 1929 by G.~D.~Birkhoff (see \cite{Grosse}): there exists an entire function $f$ whose translates are dense in the space ${\mathcal H}(\C)$ of entire functions.
In the case when $S$ is the semigroup generated by one element $s\in G$, $S$-density is equivalent to the {\it hypercyclicity} of the operator
$T_s$, which means that there exists a single vector $x$ such that the sequence $(T_s^nx)_{n\in\N}$ is dense.

If $S$ is a sub-semigroup or a subgroup of $G$, $S$-density means in other terms that the action of $S$ by translations on the space $\lpo$ admits a hypercyclic vector. On a non-weighted $L^p$-space, the translation operators are isometric and therefore cannot have hypercyclic vectors. In the weighted case
hypercyclicity can occur, depending on the weight; moreover, the conditions on weight are as a rule explicit and not difficult to calculate.
This problem was considered
by many authors. The characterisation in the discrete case $G=\Z$,  $S=\Z_+$ was obtained by H. Salas in 1995 \cite{salas}. In 1997, W. Desch, W. Schappacher, and G. Webb \cite{DSW} characterized the hypercyclicity of the $C_0$-semigroup of translations by $S=\R_+$, with
$G=\R$. The dynamics of
the translation $C_0$--semigroup indexed by a sector of the complex plane was considered by J. A. Conejero
and
A. Peris \cite{Conejero}. Ch.-Ch.~Chen \cite{Chen} characterized the hypercyclicity of a single translation operator (when $S$ is the semigroup generated by an element of $G$)  for locally compact groups.

In the present paper, we do not suppose that the set $S$ has any algebraic structure: it can be an arbitrary subset of $G$. We denote by $\bSigma$ the class of non-compact second countable locally compact groups.
For a locally compact group $G$, we show that $G\in\bSigma$ if and only if for some weight $\o$ on $G$ and for some $S\subset G$ there exist $S$-dense functions.

The main result of the paper is Theorem \ref{series-criter}, which provides a necessary and sufficient condition for the existence of $S$-dense functions in $\lpo$. To avoid technical details, we formulate here this condition in the case when $G$ is discrete:

{\bf Corollary \ref{crit-discrete}.} {\it Let $G$ be an infinite countable discrete group, and let $S\subset G$. Let $\o$ be a weight on $G$. There is an $S$-dense vector in $\lpo$ if and only if for every increasing sequence $(F_n)_{n\ge1}$ of finite subsets of $G$, there exists a sequence $(s_n)_{n\ge1}\subset S$ such that the sets $s_n^{-1}F_n$ are pairwise disjoint and, setting $s_0=e$, $F_0=\emptyset$,
$$
\sum_{n,k: n\ne k} \Big( \int_{s_ns_k^{-1}F_k} \o^p \Big)^{1/p} <\infty.
$$}

If we suppose that $S$ generates a commutative subgroup of $G$ (which is not necessarily discrete), the criterion is straightforward:

{\bf Theorem \ref{criter-abel}.} {\it
Let $G\in\bSigma$, and let $S\subset G$ generate an abelian subgroup in $G$. Let $\o$ be a weight on $G$. There is an $S$-dense vector in $\lpo$ if and only if for every compact set $F\subset G$ and any given $\delta>0$, there exist $s\in S$ and a compact $E\subset F$ such that  $\mu(F\setminus E)<\delta$ and
$$
\esssup_{sE\cup s^{-1}E} \o <\delta.
$$
}
Note that if, in addition, $G$ is discrete, compact subsets of $G$ are just finite, and $E=F$ in this condition.

Our theorems generalize the  results \cite{DSW}, \cite{salas}, \cite{Chen} mentioned above.

\section{Notations and assumptions on the group}

Recall that we consider the weighted space $L^p(G,\o)=\{f: f\o\in L^p(G)\}$ with the norm $\|f\|_{p,\o} = \big( \int|f\o|^p\big)^{1/p}$. The assumption that $\o$ is locally $p$-summable (that is $\o\in L^p(K)$ for any compact set $K\subset G$) implies that $L^p(G,\o)$ contains all measurable compactly supported bounded functions.
For a set $K\subset G$, denote by $I_K$ the characteristic function of $K$. For a measurable function $f$, denote $\|f\|_{p,K} = (\int_K |f|^p)^{1/p}$, with $\|f\|_p=\|f\|_{p,G}$.

The Haar measure is denoted by $\mu$, but we will also use the notations $|A|\equiv\mu(A)$ for $A\subset G$ and $\int_G f d\mu \equiv \int_G f(t)\,dt$. The translation $T_sf$ of a function $f\in \lpo$ will also be denoted by $\trans sf$ ($s\in G$).

For $s\in G$, one has the following expression for the norm of the translation operator:
\beq\label{norm-of-Ts}
\|T_s\| = \esssup_{t\in G} \frac{\o(st)}{\o(t)}.
\eeq
We say that a weight $\o$ is {\it $S$-admissible\/} if $\|T_s\|<+\infty$ for all $s\in S$.

Let $\mathbf\Sigma$ be the class of non-compact second countable locally compact groups. Note that every group in $\bSigma$ is $\sigma$-compact, that is, represented as a countable union of compact sets (to see this, pick an open compactly generated subgroup $H\subset G$, which is a fortiori $\sigma$-compact; $G$ being second countable, $[G:H]$ is necessarily countable, so finally $G$ is $\sigma$-compact). In this section we show that $G$ must be in $\bSigma$ to admit $S$-dense vectors for some $S\subset G$ and some weight $\o$.

In the statements \ref{lemma-not-compact}-\ref{G-sec-countable} it is supposed that $\S\subset G$ and the weight $\o$ is $S$-admissible.

Our first claim is that $G$ can't be compact. In the following lemma we prove a bit more, and it is exactly this result that we use below in the proof of Theorem \ref{series-criter}.

\blm\label{lemma-not-compact}
Let $x$ be an $S$-dense vector in $\lpo$ for some $S,p,\o$, and let $K\subset G$ be a compact set. Then for every $\e>0, \lambda\in\R$ and every compact set $L\subset G$ there is $s\in S\setminus L$ such that $\|\lambda I_K-\trans sx\|_{p,\o}<\e$.
\elm
\bpr
We can scale the Haar measure if necessary to have $\mu(K)=1$.
Increasing $L$ if necessary, we can suppose that $L=L^{-1}$.
Suppose the contrary: that is, $\|\lambda I_K-\trans sx\|_{p,\o}\ge\e$ for every $s\in S\setminus L$. There is clearly $C>0$ such that for $K_C=\{t\in K: \omega(t)>1/C\}$ we have $\mu(K_C)>3/4$.

Since $LK$ is also compact, its measure is finite. Pick $N\in \N$ with $N>2\mu(LK)$. Set $\lambda_j=\lambda + j\,\e/(2N)$, $j=0,\dots,N$ and pick $\delta>0$ such that $\delta < \e^p/(4N)^p$ and $\delta/C^p <\e$. For $j=0,\dots,N-1$, there is $s_j\in S$ such that $\|\lambda_j I_K - \trans{s_j}x\|_{p,\o}^p < \delta/(4C^p)$. Since $\|\lambda I_K-\trans{s_j}x\|_{p,\o}< \delta/(4C^p)+ |\lambda-\lambda_j|<\e/4+ \e/2<\e$, by assumption we have $s_j\in L$.

It follows that
\begin{align*}
\int_{s_j K_C} |\lambda_j-x(t)|^pdt &= \int_{K_C} |\lambda_j-x(s_j^{-1}t)|^pdt < C^p\int_{K_C} \o(t)^p|\lambda_j-x(s_j^{-1}t)|^pdt
\\&\le C^p \|\lambda_j I_K-\trans{s_j}x\|_{p,\o}^p <C^p \frac\delta{4C^p} = \frac\delta{4}.
\end{align*}
Set $K_j=\{t\in K_C: |\lambda_j-x(s_j^{-1}t)|^p\ge \delta\}$.
By trivial estimates,
$$
\frac\delta{4}>\int_{K_j} |\lambda_j-x(s_j^{-1}t)|^p dt \ge \delta\, \mu(K_j),
$$
so that
$$
\mu(K_j) < \frac14.
$$
For $t\in K'_j=s_j^{-1}(K_C\setminus K_j)$ we have the estimate $|x(t)-\lambda_j|<\delta^{1/p} < \e/(4N)$. As $|\lambda_j-\lambda_k| \ge \e/(2N)$ for $j\ne k$, it follows that the sets $K'_j$ are pairwise disjoint. Moreover, $\mu(K'_j) \ge \mu(K_C)-\mu(K_j)>1/2$ for every $j$, hence $\mu(\cup K'_j)>N/2 > \mu(LK)$ and in particular, $\cup K'_j\not\subset LK$. But by assumption, we have $K'_j\subset s_j^{-1} K\subset LK$ for every $j$. This contradiction proves the lemma.
\epr

\bcor\label{cor-G-not-compact}
If for some $p$, $S$, $\o$ there exists an $S$-dense vector in $\lpo$, then $G$ is not compact.
\ecor

\bprop\label{G-sec-countable}
If for some $p$, $S$, $\o$ there exists an $S$-dense vector in $\lpo$, then $G$ is second countable.
\eprop
\bpr
We use then the fact that $G$ is second countable if and only if $L^p(G)$ is separable \cite[Theorem 2]{devries}, and the same is valid for $\lpo$ since this space is isometrically isomorphic to $L^p(G)$.

Suppose that $\lpo$ is not separable but $x\in\lpo$ is an $S$-dense vector. There exists then an uncountable set $\{f_\a\}\subset \lpo$ such that $\|f_\a-f_\beta\|_p>3$ (this is easy to show usign the Zorn's lemma). Approximating $f_\a$ by $\trans{s_\a}x$ so that $\|f_\a-\trans{s_\a}x\|_p<1/2$, we get a family $s_\a\in S$ such that $\|\trans{s_\a}x-\trans{s_\beta}x\|_{p,\omega}^p>2$. In particular, $s_\a$ are all different.

For every $\a$,
$$
\|\trans{s_\a}x\|_{p,\omega}^p = \int \omega^p |\trans{s_\a}x|^p = \lim_{N\to\infty} \int_{\omega\le N} \omega^p |\trans{s_\a}x|^p.
$$
There exists $N_\a\in \N$ such that $\int_{\omega\ge N_\a} \omega^p |\trans{s_\a}x|^p <1/2$.

Since the set of $\a$ is uncountable, there is a single $N\in\N$ such that the set $A=\{\a: N_\a=N\}$ is uncountable. For every $\a,\beta\in A$, we have
\begin{align*}
\int_{\omega\le N} |\trans{s_\a}x-\trans{s_\beta}x|^p
& \ge N^{-p} \int_{\omega\le N} \omega^p |\trans{s_\a}x-\trans{s_\beta}x|^p
\\&\ge N^{-p} \Big( \|\trans{s_\a}x-\trans{s_\beta}x\|_{p,\omega}^p - \int_{\omega\ge N_\a} \omega^p |\trans{s_\a}x|^p - \int_{\omega\ge N_\a} \omega^p |\trans{s_\beta}x|^p \Big) > N^{-p}.
\end{align*}

But this contradicts the fact that the $G$-orbit of any function in the non-weighted $L^p(G)$ is separable \cite{tam}: pick $\gamma\in A$ and set $y=\trans{s_\gamma}x\cdot I_{\o\ge N}$, then $y\in L^p(G)$ (with norm $\le N/2^{1/p}$), and $\int_{\omega\le N} |\trans{s_\a}x-\trans{s_\beta}x|^p = \| \trans{(s_\a s_\gamma^{-1})}y-\trans{(s_\beta s_\gamma^{-1})}y\|_p^p$ for all $\a\ne \beta$ \hbox{in $A$}.
\epr

Corollary \ref{cor-G-not-compact} and Proposition \ref{G-sec-countable} prove the implication $(2)\Rightarrow(1)$ of the following

\btm\label{existence-general}
For a locally compact group $G$, the following assumptions are equivalent:
\begin{enumerate}
\item $G$ is in the class $\bSigma$ (non-compact second countable);
\item for every subset $S\subset G$ with non-compact closure there is an $S$-admissible weight $\o$ such that $\lpo$ contains $S$-dense vectors, for all (or some) $p\ge1$;
\item there is a weight $\o$ for which all left and right translations are bounded and such that $\lpo$ contains $G$-dense functions for all (or some) $p\ge1$.
\end{enumerate}
\etm
The implication $(1)\Rightarrow(3)$ is proved in the beginning of Section \ref{sec-examples} (and $(3)\Rightarrow(2)$ is obvious).

\section{Density criterion}

In the most general case, the criterion has the following form. Below we give similar conditions in several particular cases.
Recall that in every statement we assume that $T_s$ is a bounded operator on $\lpo$ for all $s\in S$.

\btm\label{series-criter}
Let $G\in\bSigma$, and let $S\subset G$. Let $\o$ be an $S$-admissible weight on $G$. There is an $S$-dense vector in $\lpo$ if and only if
for every increasing sequence $(F_n)_{n\ge1}$ of compact subsets of $G$ and any given $\delta_n>0$, there exists a sequence $(s_n)_{n\ge1}\subset S$ and compact sets $K_n\subset F_n$ such that the sets $s_n^{-1}F_n$ are pairwise disjoint, $|F_n\setminus K_n|<\delta_n$ and, setting $s_0=e$, $K_0=\emptyset$,
\beq\label{omega-series-simpl}
\sum_{n,k\ge0: n\ne k} \|\o \|^p_{p,s_ns_k^{-1}K_k} <\infty.
\eeq
\etm
\bpr
$\kern4pt\raise2.5pt\hbox{\circle{17}}\kern-5pt{\Rightarrow}$\kern5pt:\\
We set $s_0=e$ and $K_0=\emptyset$.

Decreasing $\delta_k$ if necessary, we can assume that $\sum \delta_k<\infty$. Let $y\in\lpo$ be $S$-dense. For every $k$, there is $F'_k\subset F_k$ such that $|F_k\setminus F'_k|<\delta_k/2$ and
$$
\essinf_{F'_k}\o\ge c_k, \qquad \esssup_{F'_k}\o\le c_k^{-1}
$$
with some $0<c_k<1/4$. By Lemma \ref{lemma-not-compact}, there exists $s_k\in S$ such that $\|2 I_{F'_k}-{}^{s_k}y\|_{p,\omega}^p<c_k^{2p}\,\delta_k$ and $s_k^{-1}F_k$ does not intersect $s_j^{-1}F_j$ for $j<k$ (this condition is equivalent to $s_k\notin F_k F_n^{-1} s_n$). This implies that
$$
c_k^{2p} \,\delta_k > \int_{F'_k} |{}^{s_k}y(t)-2|^p\,\o(t)^p dt \ge c_k^p \int_{F'_k} |{}^{s_k}y(t)-2|^p dt
$$
and so $\|{}^{s_k}y-2\|_{p,F'_k}^p < c_k^p\delta_k<\delta_k$.

Set $K'_k = \{t\in F'_k: |y(s_k^{-1}t)|>1\}$. Then
\beq\label{y>1}
c_k^p\delta_k > \|{}^{s_k}y-2\|_{p,F'_k}^p = \int_{K'_k} |y(s_k^{-1}t)-2|^p dt + \int_{F'_k\setminus K'_k} |y(s_k^{-1}t)-2|^p dt,
\eeq
which implies
$$
|F'_k\setminus K'_k| \le \int_{F'_k\setminus K'_k} |y(s_k^{-1}t)-2|^p dt < c_k^p\delta_k<\delta_k/4.
$$
By regularity of the Haar measure, there is a compact set $K_k\subset K'_k$ such that $|K'_k\setminus K_k|<\delta_k/4$. By estimates above, we have  $|F_k\setminus K_k|<\delta_k$.

Since the norm of $y$ is finite, we have:
$$
\sum_{k} \int_{s_k^{-1}K_k} \,\o^p \le \sum_{k} \int_{s_k^{-1}K_k} \,\o^p(t)\, |y(t)|^p dt \le \|y\|_{p,\o}^p <\infty.
$$

By the choice of $s_n$, we have
\begin{align*}
\int_{G\setminus K_n} \o^p\, |{}^{s_n}y|^p &= \int_{F'_n\setminus K_n} \o^p\, |{}^{s_n}y|^p + \int_{G\setminus F'_n} \o^p\, |{}^{s_n}y|^p
\\&\le |F'_n\setminus K_n| c_n^{-p} + \|2 I_{F'_n}-{}^{s_n}y\|_{p,\omega}^p < \delta_n + c_n^p\,\delta_n \le 2\delta_n.
\end{align*}
At the same time, $|{}^{s_n}y|^p > 1$ on $s_ns_k^{-1}K_k$. For $k\ne n$, the condition $s_n s_k^{-1}F_k \cap F_n = \emptyset$ implies that $s_ns_k^{-1}K_k\subset G\setminus K_n$; moreover, the sets $s_n s_k^{-1}F_k$, $k\ne n$, are pairwise disjoint. It follows that
\begin{align*}
2\delta_n &> \sum_{k:k\ne n} \int_{s_ns_k^{-1}K_k} \o^p\, |{}^{s_n} y|^p
> \sum_{k:k\ne n} \|\o\|^p_{p,s_n s_k^{-1}K_k}.
\end{align*}
By the choice of $\delta_n$ we have $\sum_n \delta_n<\infty$, so we get immediately \eqref{omega-series-simpl}.

$\kern4pt\raise2.5pt\hbox{\circle{17}}\kern-5pt{\Leftarrow}$\kern5pt:\\
From the assumptions it follows that $\lpo$ is separable. Choose a dense sequence $(Q_n)_{n=1}^\infty\subset \lpo$ such that every $Q_n$ is compactly supported and essentially bounded. Arrange them with repetitions in a sequence $(P_n)_{n=1}^\infty$ so that every $Q_n$ appears infinitely many times in this sequence. We will seek for $X$ ``almost'' of the form $\sum_k T_{s_k^{-1}}P_k$, with $s_n\in S$ chosen so that the series converges and $\|{}^{s_n}X-P_n\|_{p,\o}<\e_n$ for every $n$, with some $\e_n\to0$. This guarantees that $X$ is $S$-dense.

Denote $F_n=\supp P_n$ and set $p_n=\|P_n\|_\infty^p$. Pick a decreasing sequence $(\delta_n)$ such that $0<\delta_n<p_n^{-1}\, 2^{-n}$. Choose $s_n$, $K_n$ according to \eqref{omega-series-simpl}. Set
$$
a_k=\sum_{n:n\ne k} \|\o \|^p_{p,s_ns_k^{-1}K_k}.
$$
By \eqref{omega-series-simpl}, $\sum_k a_k<\infty$. Choose a subsequence $(a_{l_k})$, with $l_0=0$, such that $\sum p_k a_{l_k}<\infty$ and $P_{l_k}=P_k$ for every $k$. By the choice of $(P_k)$, this is always possible. Set $t_n = s_{l_n}$ and $E_k=K_{l_k}$, then the following series converges:
\beq\label{omega-series-K_k}
\sum_{n,k: n\ne k} p_k \|\o \|^p_{p,t_nt_k^{-1}E_k} <\infty.
\eeq
In particular, the sets $t_n^{-1}E_n$ are pairwise disjoint, and
$$
|\supp P_n\setminus E_n| = |\supp P_{l_n}\setminus K_{l_n}|<\delta_{l_n}\le \delta_n
$$
for $n\ge1$. We can still assume that $t_0=e$, $E_0=\emptyset$. Set
$$
X = \sum_k T_{t_k^{-1}}( P_k I_{E_k} ).
$$
This series converges, since (note that the support of $T_{t_k^{-1}}P_k$ is $t_k^{-1} \supp P_k$)
\begin{align*}
\|X\|^p_{p,\o}
 &= \sum_k \int |P_k(t_k t) I_{E_k}(t_kt) \o(t)|^p dt
 \le \sum_k \|P_k\|_\infty^p \int_{t_k^{-1}E_k} |\o(t)|^p dt
\\& = \sum_k p_k \|\o \|^p_{p,t_0t_k^{-1}E_k} <\infty.
\end{align*}

Now, for every $n$
\beq\label{X-Pn-1}
\| {}^{t_n}X-P_n\|_{p,\o}^p
 = \|P_n (1-I_{E_n}) \|^p_{p,\o} +  \sum_{k:k\ne n} \|{}^{t_n t_k^{-1}} ( P_k I_{E_k} ) \|^p_{p,\o},
\eeq
since again the supports of the summands in \eqref{X-Pn-1} are pairwise disjoint:
$$
\supp {}^{t_nt_k^{-1}} P_k \subset t_n t_k^{-1} E_k.
$$
One has
$$
\|P_n (1-I_{E_n}) \|^p_{p,\o} = \int_{\supp P_n\setminus E_n} |P_n|^p \le \|P_n\|_\infty^p \, |\supp P_n\!\setminus \!E_n|
 < p_n \delta_n < 2^{-n}.
$$
and for $k\ne n$,
\begin{align*}
\|{}^{t_n t_k^{-1}}(P_k I_{E_k}) \|^p_{p,\o}
 \le \|P_k\|_\infty^p \|\o\|^p_{p,t_n t_k^{-1} E_k} = p_k \|\o\|^p_{p,t_n t_k^{-1} E_k}.
\end{align*}

Denote
$$
\e_n = \sum_{n,k:k\ne n} p_k \|\o\|^p_{p,t_n t_k^{-1} E_k}.
$$
By \eqref{omega-series-K_k}, $\sum_n \e_n<\infty$, so $\e_n\to0$. Thus, we have
\beq\label{estimate-oLp-Pk}
\| {}^{t_n}X-P_n\|_{p,\o}^p \le \e_n + 2^{-n}
\eeq
with $\e_n\to0$, which proves that $X$ is $S$-dense.
\epr

\bcor\label{crit-discrete}
Let $G$ be an infinite countable discrete group, and let $S\subset G$. Let $\o$ be an $S$-admissible weight on $G$. There is an $S$-dense vector in $\lpo$ if and only if for every increasing sequence $(F_n)_{n\ge1}$ of finite subsets of $G$, there exists a sequence $(s_n)_{n\ge1}\subset S$ such that the sets $s_n^{-1}F_n$ are pairwise disjoint and, setting $s_0=e$, $F_0=\emptyset$,
\beq
\sum_{n,k\ge0: n\ne k} \|\o \|^p_{p,s_ns_k^{-1}F_k} <\infty.
\eeq
\ecor

\brem\label{remark-crit}
It is enough to check the conditions of Theorem \ref{series-criter} for an increasing sequence $(F_n)$ of compact sets such that $G=\cup_n F_n$. Indeed, for another sequence $(\tilde F_k)$ one can choose $n_k$ so that $\tilde F_k\subset F_{n_k}$, and the rest is obvious.
\erem

\section{Simplifications of conditions}

\btm\label{suff-cond-delta}
Let $G\in\bSigma$, and let $S\subset G$. Let $\o$ be an $S$-admissible weight on $G$. If for every compact set $K\subset G$ of positive measure and every $\e,\delta>0$ there exist $s\in S$ and compact $E\subset K$ such that $|K\setminus E|<\delta$ and
\beq\label{eq-suff-cond-delta}
\esssup_{\!s E\,\cup\, s^{-\!1}\!E} \; \o <\e,
\eeq
then $\lpo$ contains an $S$-dense vector.
\etm
\bpr Let us show first that under the assumptions above, $s$ can be chosen outside of a given compact set $L$. Suppose the contrary: that \eqref{eq-suff-cond-delta} does not hold for any $s\notin L$. We can assume that $\delta_1=|K|-\delta>0$. Choose a compact set $E_0$ and $s_0\in S$ so that $|K\setminus E_0|<\delta$ and \eqref{eq-suff-cond-delta} holds for $E=E_0$, $s=s_0$. In particular, $|E_0|>|K|-\delta=\delta_1$. There is a compact subset $K_1\subset E_0$ such that $|E_0\setminus K_1|<\delta_1/9$ and $\e_1=\essinf \{\o(t): t\in s_0K_1\cup s_0^{-1} K_1\} >0$.

Next choose by induction compact sets $E_n$, $K_n$ and $s_n\in S$, $n\ge1$, so that $|K_n\setminus E_n|<\delta_1/9^n$, $\esssup_{\!s_n E_n\,\cup\, s_n^{-\!1}\!E_n} \; \o <\e_n$, $|E_n\setminus K_{n+1}|<\delta_1/9^n$ and $\e_{n+1}=\essinf \{\o(t): t\in s_nK_{n+1}\cup s_n^{-1} K_{n+1}\} >0$.
For every $n\ge1$, we have
$|E_n| > |K_n|-\delta_1/9^n > |E_{n-1}|-2\delta_1/9^n$, which implies
$$
|E_n|>|E_0|- 2\delta_1\sum_{k=1}^n \frac1{9^k} > |E_0|-\delta_1/4 = 3\delta_1/4.
$$

Denote $D_n=s_n E_n\,\cup\, s_n^{-\!1}\!E_n$, $D_n^- = \cup_{k<n} D_k\cap D_n$ and $D_n^+ = \cup_{k>n} D_k\cap D_n$. Since $D_k\cap D_n$ is contained (except probably for a zero measure subset) in $(\{s_n\}\cup \{s_n^{-1}\})(E_n\setminus K_{n+1})$ for all $k>n$, we have $|D_n^+| \le 2|E_n\setminus K_{n+1}| < 2\delta_1/9^n$. As a consequence,
$$
|D_n^-| \le \sum_{k<n} |D_k^+| < 2\delta_1\sum_{k=1}^\infty \frac1{9^k} = \frac{\delta_1}{8},
$$
and for every $n$,
\begin{align*}
|\cup_{k\le n} D_k | &= |\cup_{k<n} D_k \sqcup (D_n\setminus D_n^-) |
 \\&\ge |\cup_{k<n} D_k| +|D_n|-|D_n^-| \ge |\cup_{k<n} D_k|+\frac{3\delta_1}{4}-\frac{\delta_1}{8} > |\cup_{k<n} D_k|+\frac{\delta_1}{2},
\end{align*}
whence it follows that
$$
| LE_0\cup L^{-1}E_0 | \ge |\cup_{k=1}^\infty D_k |\ge \sum_{k=1}^\infty \frac{\delta_1}{2} = +\infty,
$$
which is impossible since $LE\cup L^{-1}E$ is compact. This contradiction shows that $s$ can be chosen outside $L$.

Now let us show that the assumptions of the theorem imply the assumptions of Theorem \ref{series-criter}.

Let $(F_n)$ be an increasing sequence of compact sets and $(\delta_n)$ a sequence of positive numbers. We can assume that $\delta_1<1/2$ and $\delta_{n+1}<\delta_n/2$ for all $n$, in particular $\delta_n<2^{-n}$.
Choose firstly sequences $E'_n\subset G$, $s_n\in S$ by induction. Set $E'_0=\emptyset$, $s_0=e$; for $n\ge1$, set
$$
K_n = \Big(\bigcup_{k<n} s_kF_n\Big)\cup \Big(\bigcup_{k<n} F_k \cup s_k^{-1}F_k\Big).
$$
Set $C_n = \max_{k<n} \|T_{s_k}\|$ (automatically $C_n\ge1$)
and choose $s_n\in S$, $E'_n\subset K_n$ so that $E'_n$ is compact, $|K_n\setminus E'_n|<\delta_n/2$ and $\|\o \|_{\infty,s_nE'_n\cup s_n^{-1}E'_n} <(2^{n}C_n|E'_n|^{1/p})^{-1}$. We can choose $s_n$ so that $s_n^{-1}F_n$ does not intersect $s_k^{-1}F_k$ for $k<n$ (that is, $s_n\notin F_n F_k^{-1} s_k$).

Set now $E_0=\emptyset$ and
\beq\label{def-En}
E_n = F_n \cap E_n'\cap s_n(\cap_{k> n} E'_k)
\eeq
for $n\ge1$.
Then, since $F_n\subset K_n$ and for $k>n$ \ $s_n^{-1}F_n\subset K_k$, that is $F_n\subset s_nK_k$, we have
\begin{align*}
|F_n\setminus E_n| & \le |K_n\setminus E'_n| + |\cup_{k>n} s_nK_k\setminus s_nE'_k|
 = \sum_{k\ge n} |K_k\setminus E'_k|
\\& <\sum_{k\ge n} \frac{\delta_k}2
<\sum_{k\ge n} \frac{\delta_n}{2^{k-n+1}} = \delta_n.
\end{align*}
The sets $s_n^{-1}F_n$ are pairwise disjoint; for $k<n$ we have $E_k\subset s_kE'_n$ by \eqref{def-En}, which implies $s_k^{-1}E_k\subset s_k^{-1}F_k\subset E'_n$. Now
\begin{align*}
\sum_{n,k: n\ne k} \|\o \|^p_{p,s_ns_k^{-1}E_k}
&= \sum_n \bigg( \sum_{k<n} \|\o \|^p_{p,s_ns_k^{-1}E_k} + \sum_{k>n} \|\o \|^p_{p,s_ns_k^{-1}E_k}\bigg)
\\&\le \sum_n \bigg( \|\o \|^p_{p,s_nE'_n} + \sum_{k>n} \|T_{s_n}\|^p \|\o \|^p_{p,s_k^{-1}E_k}\bigg)
\\&\le \sum_n \bigg( \|\o \|^p_{\infty,s_nE'_n} |E'_n| + \sum_{k>n} C_k^p  \|\o \|^p_{\infty,s_k^{-1}E'_k} |E'_k| \bigg)
\\&= \sum_n \bigg( 2^{-n} + \sum_{k>n} 2^{-k}\bigg) <\infty.
\end{align*}
Finally, the assumptions of Theorem \ref{series-criter} are satisfied with $E_n$ in place of $K_n$.
\epr

\bcor\label{suff-cond}
Let $G\in\bSigma$, and let $S\subset G$. Let $\o$ be an $S$-admissible weight on $G$. If for every compact set $K\subset G$ of positive measure
\beq\label{inf-esssup-omega}
\inf_{ s\in S} \;\;\Big(\esssup_{\!s K\,\cup\, s^{-\!1}\!K} \; \o \Big) = 0,
\eeq
then $\lpo$ contains an $S$-dense vector.
\ecor

In general, this condition is not necessary, see Example \ref{ex-hyper-but-not-suff-cond}. But if $S$ generates an abelian subgroup, this is the case:

\btm\label{criter-abel}
Suppose that $G\in\bSigma$, and let $S\subset G$ generate an abelian subgroup. Let $\o$ be an $S$-admissible weight on $G$. The following conditions are equivalent:
\begin{enumerate}
\item There is an $S$-dense vector in $\lpo$ for every $p$, $1\le p<\infty$;
\item There is an $S$-dense vector in $\lpo$ for some $p$, $1\le p<\infty$;
\item For every compact set $F\subset G$ and any given $\delta>0$, there exist $s\in S$ and a compact $E\subset F$ such that  $|F\setminus E|<\delta$ and
\beq\label{eq-criter-abel}
\esssup_{sE\cup s^{-1}E} \o <\delta;
\eeq
\item For some $p$, $1\le p<\infty$, the following condition holds: For every compact set $F\subset G$ and any given $\delta>0$, there exist $s\in S$ and a compact $E\subset F$ such that  $|F\setminus E|<\delta$ and
\beq\label{eq-criter-abel-lp}
\|\o\|_{p,sE\cup s^{-1}E} <\delta.
\eeq
\item For every $1\le p<\infty$, the assumptions of (4) above hold.
\end{enumerate}
\etm
\bpr
(4)$\Rightarrow$(3) is an easy exercise, (3)$\Rightarrow$(1) was shown in Theorem \ref{suff-cond-delta} without the commutativity assumption; (1)$\Rightarrow$(2) and $(5)\Rightarrow(4)$ are obvious.

It remains to prove (2)$\Rightarrow$(5).\\
If there is an $S$-dense vector in $\lpo$, then for $F_n=F$ and $\delta_n=2^{-n}$, $n\ge1$, there exist $s_n\in S$, $E_n\subset F$ with $s_0=e$, $E_0=\emptyset$ such that $|F\setminus E_n|<\delta_n$ and
$$
\sum_{n,k: n\ne k} \|\o \|^p_{p,s_ns_k^{-1}E_k} <\infty.
$$
It follows (by considering the subseries with $n=0$) that $\|\o \|^p_{p,s_k^{-1}E_k}\to0$.

Let $k$ be such that $\delta_k<\delta/2$. Set $C=\|T_{s_k}\|$. For $n>k$, we have
$$
\|\o\|^p_{p,s_nE_k} = \int_{s_nE_k} |\o|^p = \int_{s_ks_ns_k^{-1}E_k} |\o|^p
\le C \int_{s_ns_k^{-1}E_k} |\o|^p \to0, \quad n\to\infty.
$$
Pick now $n>k$ such that $\|\o \|_{p,s_n^{-1}E_n}<\delta/2$ and $\|\o\|_{p,s_nE_k}<\delta/2$. Set $E=E_k\cap E_n$, $s=s_n$, then $|F\setminus E|<\delta_k+\delta_n<\delta$ and $\|\o \|_{p,sE\cup s^{-1}E} \le \|\o\|^p_{p,s_nE_k} + \|\o\|^p_{p,s_n^{-1}E_n}<\delta$.
\epr

The following case is quite particular and is meaningful especially in the abelian case, see Corollary \ref{abelian}.

\brem\label{equiv-to-continuous}
If $\o$ is locally summable and
\beq\label{L}
L(g)=\sup_{t\in G} \frac{\o(gt)}{\o(t)}<+\infty \qquad
\eeq
for all $g\in G$, then $\o$ is equivalent to a continuous weight \cite[Satz 2.7]{feicht}. Since \eqref{L} holds or does not hold simlutaneously for $\o^p$ for all $1\le p<\infty$, we have the same conclusion if $\o$ is locally $p$-summable and satisfies \eqref{L}. Hence in the following statement we can assume, without loss of generality, that $\omega$ is continuous.
\erem

\bprop\label{all-equiv-to-inf}
Let $G\in\bSigma$, and suppose that the weight $\o$ is continuous and such that for all $g\in G$,
\beq\label{R}
L(g)=\|T_g\|=\sup_{t\in G} \frac{\o(gt)}{\o(t)}<+\infty, \qquad
R(g)=\sup_{t\in G} \frac{\o(tg)}{\o(t)}<+\infty.
\eeq
Then for any $S\subset G$, the following conditions are equivalent:
\begin{enumerate}
\item There are $S$-dense vectors in $\lpo$.
\item
\beq\label{inf-sup-omega}
\inf_{ s\in S} \;\;\sup_{\!s K\,\cup\, s^{-\!1}\!K} \; \o  = 0 \text{ \ for every compact set $K\subset G$}
\eeq
\item
\beq\label{inf-omega}
\inf_{s\in S} \max\big( \o(s), \o(s^{-1}) \big) = 0
\eeq
\end{enumerate}
\eprop
\bpr
Obviously \eqref{inf-sup-omega}$\Rightarrow$\eqref{inf-omega}. Let us prove the inverse implication.
Let $K$ be a compact set; we can assume that $K=K^{-1}$. It is proved in a rather general situation in \cite[Proposition 1.16]{edwards} that $L$ and $R$ are locally bounded (provided they are finite). There exists thus a constant $C$ such that $L|_K \le C$ and $R|_K\le C$.

Now for every $s\in G$,
\begin{align*}
\sup_{t\in sK} \o(t) &
\le \sup_{t\in K} R(t)\o(s) \le C \o(s).
\end{align*}
In the same way we get the inequality $\sup_{s^{-1}K} \o \le C\o(s^{-1})$, thus \eqref{inf-sup-omega} and \eqref{inf-omega} are equivalent.

Since they both imply \eqref{inf-esssup-omega}, they imply also the assumption \eqref{omega-series-simpl} of Theorem \ref{series-criter} and the existence of $S$-dense vectors in $\lpo$.

Suppose now that there exist $S$-dense vectors in $\lpo$. Pick a compact set $F$ of positive measure and set $F_n=F$ for all $n\in\N$. Choose $s_n\in S$ and $K_n\subset F$ so that \eqref{omega-series-simpl} holds. By setting $\delta_n=|F|/4$, we can guarantee that $|K_n\cap K_1|>|F|/2$ for all $n\ge 1$. From \eqref{omega-series-simpl} it follows, in particular, that
$$
\lim_{n\to\infty} \|\o \|^p_{p,s_ns_1^{-1}K_1} = 0.
$$
Let $C_1>0$ be such that $L(s_1^{-1})\le C_1$, $R|_{K_1^{-1}}\le C_1$ and
$R|_{K_1^{-1}s_1}\le C_1$. We have
$$
\o(s_n) = \inf_{t\in s_1^{-1} K_1} \o(s_n tt^{-1}) \le \inf_{t\in s_1^{-1}K_1} R(t^{-1})\o(s_n t) \le C_1 \|\o \|_{p,s_ns_1^{-1}K_1} |K_1|^{-1}
$$
which implies $\o(s_n)\to0$, $n\to\infty$.
At the same time,
\begin{align*}
\o(s_n^{-1}) &= \inf_{t\in K_1\cap K_n} \o(s_1^{-1} s_1 s_n^{-1} tt^{-1}) \le \inf_{t\in K_1\cap K_n} L(s_1^{-1}) \o(s_1 s_n^{-1} t) R(t^{-1})
 \\&\le C_1^2 \inf_{t\in K_1\cap K_n} \o(s_1 s_n^{-1} t) \le C_1^2 \|\o \|_{p,s_1s_n^{-1}K_n} |K_1\cap K_n|^{-1}
  < 2 C_1^2 |F|^{-1} \|\o \|_{p,s_1s_n^{-1}K_n} .
\end{align*}
By \eqref{omega-series-simpl}, $\o(s_n^{-1})\to0$, $n\to\infty$, and this implies \eqref{inf-omega}.
\epr

In the case when $G$ is abelian our conditions take an especially simple form:
\bcor\label{abelian}
Let $G\in\bSigma$ be abelian, $1\le p<\infty$, and let $\o$ be a continuous $G$-admissible weight.
For any given $S\subset G$, there are $S$-dense vectors in $\lpo$ if and only if
\beq\label{inf-omega-on-G}
\inf_{s\in S} \max\big( \o(s), \o(s^{-1}) \big) = 0.
\eeq
\ecor

\section{Examples and counter-examples}\label{sec-examples}

First we show that every group $G\in \bSigma$ admits a weight with $G$-dense vectors; this completes the proof of Theorem \ref{existence-general}.

\bpr
Let $U=U^{-1}$ be a neighbourhood of identity with compact closure and let $G_1$ be the subgroup generated by $U$. Pick $g_n\in G$, $n\in\N$, such that $G = \cup_n (g_n G_1)$; this is possible since $G$ is $\sigma$-compact. Set $U_1=U$ and for $n>1$ set by induction
$$
U'_n = \{g_n\} \cup \big( \mathop{\cup}\limits_{k=1}^{n-1} U_k U_{n-k}\big), \qquad U_n = U'_n \cup (U'_n)^{-1}.
$$
In particular, we have $U^n \subset U_n$ and $G = \cup_n U_n$. By construction, $U_j U_k\subset U_{j+k}$ for all $j,k$.

Every $U_n$ has compact closure, so we can choose $s_n\in S$, $n\in\N$, so that $E_n=U_{n+1}s_n U_{n+1} \cup U_{n+1}s_n^{-1}U_{n+1}$ are pairwise disjoint and $U_{n+1}s_n U_{n+1} \cap U_{n+1}s_n^{-1}U_{n+1}=\emptyset$. Denote $V_{n,k} = U_k s_n U_k \cup U_k s_n^{-1} U_k$ for $1\le k\le n$, $V_{n,k}=G$ for
$k>n$ and $V_{n,0}=\{s_n\}\cup\{s_n^{-1}\}$.
We have $U_j\,(V_{n,k}\setminus V_{n,k-1}) U_j \subset V_{n,k+j}\setminus V_{n,k-j-1}$ for all $n,k$.
Set
$$
\omega(s) = \begin{cases}
1, & s\notin \bigcup_n V_{n,n},\\
2^{-n+k}, & s\in V_{n,k} \setminus V_{n,k-1} \ \ (k\le n).
\end{cases}
$$
We have thus $\o(s_n)=\o(s_n^{-1})=2^{-n}$ and $\o(t)\le1$ for all $t\in G$.

Let $t$ be in $U_j$. If $s\notin \cup_n V_{n,n}$, then $st\notin \cup_{n;\, k\le n-j-1} V_{n,k}$ and $\o(st)\ge 2^{-j}$. This implies
$$
\frac1{2^j} \le \frac{\o(st)}{\o(s)} \le 1.
$$
Similarly $2^{-j} \le \o(ts)/\o(s) \le 1$.
If $s\in V_{n,k} \setminus V_{n,k-1}$ with $k\le n$, then $\o(s)=2^{k-n}$ and $st\in V_{n,k+j}\setminus V_{n,k-j-1}$, thus
$$
\frac1{2^j} = \frac{2^{k-j-n}}{2^{k-n}} \le \frac{\o(st)}{\o(s)} \le \frac{2^{k+j-n}}{2^{k-n}} = 2^j.
$$
Similarly $2^{-j} \le \o(ts)/\o(s) \le 2^j$, and we see that these inequalities hold for all $s\in G$, $t\in U$.

It follows that in the notations of \eqref{R}, $L(t)<+\infty$, $R(t)<+\infty$ for all $t\in U_j$, and since $G=\cup_j U_j$, this is true for all $t\in G$. By construction, the condition \eqref{inf-omega} holds, so we can conclude that $\lpo$ contains an $S$-dense vector.
\epr

\def\F{{\mathbb F}}
Let $\F_2$ be the free discrete group of two generators $a$ and $b$.

\bex
There exists an $\F_2$-admissible weight on $\F_2$ which satisfies the condition \eqref{inf-omega} but does not admit $\F_2$-dense vectors.
\eex
Set $U=\{e,a,b,a^{-1}, b^{-1}\}$, then $G=\cup_n U^n$. Set $n_k=2^{2^k}$, $k=1,2,\dots$. Define the weight $\o$ as follows: $\o(a^{\pm n_k})=2^{-k}$,
$\o(t)=2^{j-k}$ if $t\in (U^j \setminus U^{j-1})\,a^{\pm n_k}$, $1\le j\le k$, \hbox{and 1} elsewhere. In particular, $\o|_U=1$.

For $\a\in\{a,b,a^{-1},b^{-1}\}$, denote by $F_\a$ the set of reduced words in $\F_2$ which end by $\a$. Note that $U^k a^{\pm n_k}\subset F_a\cup F_{a^{-1}}$ for all $k$ (since every word in $U^k$ has length at most $k<n_k$). It follows that $\o\equiv 1$ on $F_b\cup F_{b^{-1}}$.

Obviously $\o(sg)\le 2\o(g)$ for $s\in U$ and for all $g\in \F_2$. It follows that $\o$ is of moderate growth in the sense of Edwards \cite{edwards}, that is, $\|T_g\|<\infty$ for all $g\in \F_2$. Also, $\inf_g \max( \o(g),\o(g^{-1})) \le \inf_k \max(\o(a^{n_k}), \o(a^{-n_k})) = 0$.

Suppose now that there is a $\F_2$-dense vector $x\in \lpo$. For every $0<\e<1/2$, there exists $g\in\F_2$ such that $\|I_U-\trans gx\|_{p,\o}<\e$. This implies that
$$
|x(g^{-1}t)-1|\o(t) <\e
$$
for $t=e,b$. It follows that $|x(t)|>1/2$ for $t=g^{-1}$ and $t=g^{-1}b$.

By Lemma \ref{lemma-not-compact}, for the same $\e$ there is $h\in \F_2\setminus (Ug\cup Ub^{-1}g)$ such that ${\|{}^hx-I_U\|_{p,\o}<\e}$. This implies that $|x(h^{-1}t)|\,\o(t)<\e$ for $t\notin U$. In particular, setting $t=hg^{-1}$ we get:
$|x(h^{-1}hg^{-1})|\,\o(hg^{-1}) = |x(g^{-1})|\,\o(hg^{-1})<\e$ which implies $\o(hg^{-1})<2\e<1$; setting $t=hg^{-1}b$, we infer that $\o(hg^{-1}b)<1$. By construction of $\o$, this implies that $hg^{-1},hg^{-1}b\notin F_b\cup F_{b^{-1}}$, which is impossible. This contradiction proves that there are no $\F_2$-dense vectors.

\bex\label{ex-hyper-but-not-suff-cond}
There exists an $S$-admissible weight on $G=\F_2$ which admits $S$-dense vectors in $\lpo$ for some semigroup $S$ but does not satisfy the sufficient condition \eqref{inf-esssup-omega}.
\eex
Set $U=\{e,a,b,a^{-1}, b^{-1}\}$, then $G=\cup_n U^n$. Set $n_k=2^{2^k}$, $k=1,2,\dots$.

For $k\in\N$, set $s_k = a^{n_k}b$ and let $S$ be the semigroup generated by $\{s_k:k\in\N\}$. Set $V_{l,k} = s_l s_k^{-1} U^k$.
Every element of $V_{l,k}$, $l\ne k$, has the form $s_l s_k^{-1} t = a^{n_l-n_k} t$ with some $t\in U^k$; in the reduced form it equals $a^p r$ with $|p| \ge |n_l-n_k|-k >0$ and $r\in U^{2k}$ not starting with $a^{\pm1}$.

Let us show that the sets $SV_{l,k}$, $k\ne l$, are pairwise disjoint. For $s=a^{n_{j_1}}b\dots a^{n_{j_m}}b$, $s'=a^{n_{j'_1}}b\dots a^{n_{j'_{m'}}}b$ and $a^p r\in V_{l,k}$, $a^{p'}r'\in V_{l',k'}$ (the latter two reduced as above) the equality $sa^pr = s'a^{p'}r'$ would mean
$$
a^{n_{j_1}}b\dots a^{n_{j_m}}b a^pr = a^{n_{j'_1}}b\dots a^{n_{j'_{m'}}}b a^{p'}r',
$$
and clearly both sides are in their reduced form. We can suppose that $m\le m'$. It follows that $j_i=j'_i$ for $i=1,\dots,m$. The equality $p=n_{j_{m+1}}$ is in fact impossible: it would imply $p>0$, then $l>k$; but in this case $p\ge n_l-n_k-k > n_l/2 > n_{l-1}$ and $p\le n_l-n_k+k<n_l$, so $p$ cannot be equal to $n_\nu$ for any $\nu\in\N$. We conclude that $m=m'$, $s=s'$ and $p=p'$.

As discussed above, $p=n_l-n_k+\xi$ with $|\xi|\le k$, and $p'=n_{l'}-n_{k'}+\xi'$ with $|\xi'|\le k'$. Either $p=p'>0$ and then $k>l$, $k'>l'$, or $p=p'<0$ and $k<l$, $k'<l'$. In the first case we can suppose that $l\ge l'$. We have $n_l = n_{l'}+n_k-n_{k'}+\xi'-\xi$; by the choice of $(n_k)$, an equality is possible only if $l=l'$. But then $n_k=n_{k'}+\xi-\xi'$, and for the same reason it follows that $k=k'$. In the second case one arrives similarly at the same conclusion. Finally, $l=l'$ and $k=k'$, which shows that $V_{l,k}$ are pairwise disjoint indeed.

Set now
$$
\o(t) = \begin{cases} 8^{-l-k}, &t\in S V_{l,k};\\
1, & \text{otherwise}
\end{cases}
$$
By the reasoning above, $\o$ is well defined.
We have
$$
\sum_{n,k: n\ne k} \|\o \|^p_{p,s_ns_k^{-1}U^k} \le \sum_{n,k: n\ne k} 8^{-n-k}(4^{k}+1) <\infty,
$$
and this implies that the assumptions of Theorem \ref{series-criter} hold.

It remains to show that $\|T_{s_n}\| < \infty$ for all $n$.
If $t\notin \cup_{l,k} S V_{l,k}$, then $\o(t)=1$ and $\o(s_nt)/\o(t) \le 1$.
If $t = s s_l s_k^{-1} u$ with $s\in S$ and $u\in U^k$, then $\o(s_nt) = \o(s_ns s_l s_k^{-1} u) = \o(t)$. It follows that $\|T_{s_n}\|\le1$.

\bex
In general, one cannot have $G_n=F_n$ in Theorem \ref{series-criter}, so that
\beq\label{inf-lpnorm-omega-no-sequence}
\sum_{n,k: n\ne k} \|\o \|^p_{p,s_ns_k^{-1}F_k} <\infty.
\eeq
This is shown below: there exists a weight $\omega$ on $G=\R^2$ which is $S$-admissible and such that the condition \eqref{eq-criter-abel} is satisfied for $S=\R\times\{0\}$, i.e. $\lpo$ contains an $S$-dense vector, but the strengthened condition \eqref{inf-lpnorm-omega-no-sequence} does not hold.
\eex
Denote $\ell_0=0$ and $\ell_n = \sum_{k=1}^n k =n(n+1)/2$ for $n\ge1$.
Define $\o$ as follows: $\o(-x,y)=\o(x,y)$, and for $x\ge0$
$$
\o(x,y) = \begin{cases}
    2^{-n}, &x\in[\ell_n,\ell_{n+1}), \; y\notin(1-2^{1-n},1), \; n\in\Z\\
    2^{n}, &x\in[\ell_n,\ell_{n+1}), \; y\in[1-2^{-n},1), \; n\in\Z\\
    2^{n-2k}, &x\in[\ell_n+k,\ell_n+k+1), \; y\in(1-2^{1-n},1-2^{-n}), \; n\in\Z, 0\le k<n
\end{cases}
$$
This is illustrated as follows:
$$
\begin{tabular}{r|c|c|c|c|c|c|c|c|c}
\multirow{2}{*}{$\scriptstyle y=1$}&\vrule height20pt width0pt$2^{1-n}$&\multicolumn{5}{c|}{$2^{-n}$}&\multicolumn{3}{c}{$2^{-n-1}$}\\
\cline{2-10}
\multirow{2}{*}{$\scriptstyle y=1-2^{-n-1}$}& \multirow{3}{*}{$2^{n-1}$}&\multicolumn{5}{c|}{\multirow{2}{*}{$2^{n}$}}&\multicolumn{3}{c}{$2^{n+1}$}\\
\cline{8-10}
\multirow{2}{*}{$\scriptstyle y=1-2^{-n}$}& &\multicolumn{5}{c|}{ }&$2^{n+1}$&$2^{n-1}$&$\cdots$\\
\cline{3-10}
\multirow{2}{*}{$\scriptstyle y=1-2^{1-n}$}& \vrule width 40pt height0pt depth0pt&\vrule height20pt width0pt $2^{n}$&$2^{n-2}$&$\cdots$&$2^{2-n}$&$2^{-n}$&\multicolumn{3}{c}{\multirow{2}{*}{$2^{-n-1}$}}\\
\cline{2-7}
\vrule height20pt width0pt depth 10pt& $\cdots$&\multicolumn{5}{c|}{$2^{-n}$}\\
\multicolumn{2}{c}{$\scriptstyle \quad \ x=\ell_{n-1}$}&\multicolumn{2}{c}{$\scriptstyle \kern-45pt x=\ell_{n}$}&\multicolumn{6}{c}{$\scriptstyle \quad \ x=\ell_{n+1}$}\\
\end{tabular}
$$
For all $t\in [0,1)$,
$$
\sup_{x,y\in \R^2} \frac{\o(x+t,y)}{\o(x,y)} = \sup_{n\in\Z,y\in[0,1]} \frac{\o(n+1,y)}{\o(n,y)};
$$
it is clear from the picture above that this ratio is always bounded by 4.
It follows that $\|T_{(t,0)}\|=4$ for all $t\in[0,1)$, so that $S=\R\times\{0\}$ acts continuously on $L^p(\R^2)$.

Setting $F_n = [n,n+1]\times [0,1]$, we have $\int_{F_n} \o^p \ge1$ for all $n$, so clearly the condition \eqref{inf-lpnorm-omega-no-sequence} does not hold.

To check \eqref{eq-criter-abel}, it is sufficient to consider $F=[-n,n]\times[-n,n]$.
Choose $t\in\N$ so that
$$
t\ge2, \quad 2^{-t}<\delta, \quad n2^{3-t}<\delta.
$$
 Set $E = F \setminus ([-n,n]\times(1-2^{2-t},1))$.
Then $|F\setminus E| = n2^{3-t}<\delta$.

Set $s=(2\ell_{t},0)$.
Since $\ell_t>t$,
$$
s+E \subset [\ell_t,+\infty)\times(\R\setminus (1-2^{2-t},1))
$$
and
$$
-s+E \subset (-\infty,-\ell_{t}]\times(\R\setminus [1-2^{1-t},1)).
$$
This implies
$$
\sup_{(s+E)\cup(-s+E)} \o^p \le 2^{-{t}} <\delta,
$$
thus \eqref{eq-criter-abel} holds.

If $S$ generates $G$ and $G$ is abelian, then such an example is impossible, see Proposition \ref{all-equiv-to-inf}.

In \cite{Conejero}, it was shown on an example that a semigroup indexed by a complex sector can be hypercyclic without any single operator being hypercyclic. We provide a simple example illustrating the same situation:
\bex
There is a $\Z^2$-admissible weight on $\Z^2$ with $\Z^2$-dense vectors such that no single operator $T_{(m,n)}$ is hypercyclic. The weight is defined as follows: set
$\o(l,m)=2^{k-n}$ if $\max(|l-n|,|m-2^n|)=k$ or $\max(|l+n|,|m+2^n|)=k$ for some $n\in\N$ and $0\le k\le n$, and $\o(l,m)=1$ otherwise.

It is easy to check that all the translations $T_{(l,m)}$, $(l,m)\in\Z^2$, are bounded. By Corollary \ref{abelian}, the action of $\Z^2$ on $L^p(\Z^2)$ for every $1\le p<\infty$ is hypercyclic. But for fixed $l,m$, there are only finitely many points on the line $(l,m)\Z$ in which $\o$ is different from 1; thus, the action of $S=\{T_{(l,m)}^n\}_{n\in\N}$ is not hypercyclic, by Proposition \ref{all-equiv-to-inf}.
\eex


\begin{thebibliography}{99}
\bibitem{Chen}
Chen, Ch.-Ch. Hypercyclic weighted translations generated by non-torsion elements. {\it Arch. Math. (Basel)} {\bf 101} (2013), no.~2, 135--141.
\bibitem{Conejero}
Conejero, J.~A.; Peris, A. Hypercyclic translation $C_0$-semigroups on complex sectors. {\it Discrete Contin. Dyn. Syst.} {\bf 25} (2009), no.~4, 1195--1208.
\bibitem{conejero2}
Conejero, J.~A.; M\"uller, V.; Peris, A. Hypercyclic behaviour of operators in a hypercyclic $C_0$-semigroup. {\it J. Funct. Anal.} {\bf 244} (2007), no.~1, 342--348.
\bibitem{devries} de Vries, J. The local weight of an effective locally compact transformation group and the dimension of $L^2(G)$.
{\it Colloq. Math.} 39 (1978), no.~2, 319--323.
\bibitem{DSW}
Desch, W.; Schappacher, W.; Webb, G.~F. Hypercyclic and chaotic semigroups of linear operators. {\it Ergodic Theory Dyn. Syst.} {\bf 17} no.~4 (1997), 793--819.
\bibitem{edwards}
Edwards R.~E. The stability of weighted Lebesgue spaces.
{\it Trans. Amer. Math. Soc.}
{\bf 93} (1959), 369--394.
\bibitem{feicht}
Feichtinger H. G.
Gewichtsfunktionen auf lokalkompakten Gruppen,
{\it Sitzber. \"Osterr. Akad. Wiss. Abt. II,}
{\bf 188}, no.~8--10 (1979), 451--471.
\bibitem{Grosse}
Grosse-Erdmann, K.-G. Universal families and hypercyclic operators. {\it Bull. Amer. Math. Soc. (N.S.)} {\bf 36} (1999), no.~3, 345--381.
\bibitem{HR}
{Hewitt, E.; Ross, K.A.}
{\it Abstract harmonic analysis I, II.}
{Springer--Verlag, 3rd printing, 1997.}
\bibitem{salas}
{Salas, H.,}
{\it Hypercyclic weighted shifts,}
{\it Trans. Amer. Math. Soc.}
{\bf 347} (1995), 993--1004.
\bibitem{tam}
Tam K.~W. On measures with separable orbit.
{\it Proc. Amer. Math. Soc.} {\bf 23} (1969), 409--411.
\end{thebibliography}
\end{document}